\def \bc {\begin{center}}
\def \ec {\end{center}}

\def \bfr {\begin{flushright}}
\def \efr {\end{flushright}}

\def \v {\vskip}
\def \ii {\'\i}

\def \ba {\begin{array}}
\def \ea {\end{array}}

\def \bea {\begin{eqnarray}}
\def \eea {\end{eqnarray}}

\def \be {\begin{equation}}
\def \ee {\end{equation}}

\def \e {\hbox{e}}

\documentstyle[12pt,a4]{article}
\begin{document}
%

\begin{center}

{\bf CAN OUR NUMBER SYSTEM BE IMPROVED?}

\vskip1cm
 {\it Miguel Navarro}\footnote{e-mail: mnavarro@ugr.es}
\end{center}
\v4mm

\noindent Instituto Carlos I de F\ii sica Te\'orica y
Computacional. Universidad de Granada, Spain. \footnote{On leave
of absence.
}

\v5mm \centerline{\bf Abstract} \v2mm \footnotesize

\noindent Our number system is a magnificent tool. But it is far
from perfect. Can it be improved? In this paper some possibilities
are discussed, including the use of a different base or directed
(negative as well as positive) numerals. We also put forward some
suggestions for further research.

 \v3mm 
\noindent Keywords: General Mathematics,
Number Theory, Mathematics Education. \normalsize
\vskip3cm \setcounter{page}{1}
\section{Introduction}

We are so accustomed to our number system that we hardly notice
it. It is not possible however to lay attention on this tool and
not to feel captivated by its elegance, power and usefulness. If
there is art in Mathematics there is no doubt that our number
system is one of its best pieces. But if our number system has art
in it, it is immensely useful and practical too.

The system we use to write number is usually called the 'decimal'
system (Fig. 1). But the fact that it uses ten symbols is not its
main or, even less so, defining characteristic. A more fundamental
feature is the use it makes of the place value of the numerals.
This place-value notation, which is easier to use than to explain
- but we will come back to it later - implies, for instance, that
12 is not the same as 21.

Our number system is surprisingly modern. Its elements appeared in
different historic moments in China and India and propagate to the
West through the muslim world. But the system was not commonly
used (by the educated elites) in the West until the dawn of the
Renaissance. There must be no doubt that the introduction of this
number system, so convenient, have something to do with the quick
developments in Mathematics, the sciences and technology from the
Renaissance onwards. It is difficult to imagine that the calculus,
for instance, would have been developed had Mathematicians have to
wrestle with a system so awkward as the Roman.

\bc
\begin{tabular}{c c c c c c c c c c c c c c c c }
  \hline &  ${}^|$&${}^|$ &  ${}^|$&${}^|$  & ${}^|$ & ${}^|$ & ${}^|$ &${}^|$
  &  ${}^|$ &${}^|$  &${}^|$  &${}^|$  &${}^|$  & ${}^|$&${}^|$  \\
  & $-4$ & $-3$
& $-2$ & $-1$
& $\>0$ & $\>1$ & $\>2$ & $\>3$ & $\>$4  & $\>5$ & $\>6$&$\>7$ &$\>8$&$\>9$ & $10$\\
& & & & & & & & & & & & & & &
\\
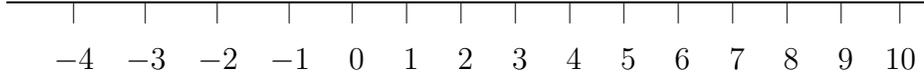
\multicolumn{14}{l}{\footnotesize{Figure 1: The number line in our number system }}\\
\end{tabular}
\ec

It is surprising that the place value principle were introduced so
late to write numbers considering that this same principle plays a
principal role in language too. In fact, it is thanks to this
principle that with a handful of letters and another symbols we
are able to express an immense complexity of meanings. Thanks to
place value the word 'top', for instance, can mean something
completely different from the word 'pot' [1]. This principle of
order also plays a role (minor perhaps in most languages) at a
grammatical level. A different word order may produce a different
meaning for the phrase. For example, 'You are tired' is an
affirmation whereas 'Are you tired?' is a question [2].

\section{Is base 6 better?}

Our number system is so useful, so convenient, so practical for
the function it is meant to fulfil that it is not possible to
think about it and not to feel fascinated.

Not everything is perfect about this tool however. To start with
children needs (many) years to learn the rudiments of the times
table. And most of them reach adulthood without completely
mastering it. I myself, in spite of having spent decades doing
Mathematics at all levels, need to use all sorts of tricks (or a
calculator) to find, for instance, what  $7\times 8$ is (exactly).

This difficulty in being mastered shed a shadow of doubt on the
perfection of a tool which the modern man needs to use almost in a
daily basis.

And this difficulty is not all. Our number system has several
other problems and limitations. The fact is, for instance, that
with our system we can write a 'few' numbers only. For instance,
no irrational number can be written with this system. To denote
irrational quantities we have to use other symbols, e.g. $\pi$,
$\e$, $\sqrt{2}$, and so on.

On the other hand, there exist the general impression that decimal
numbers are more or less the same as fractions, i.e. the rational
numbers. The truth however is that most of the common fractions,
i.e. the fractions we use on a daily basis - such as $\frac13$,
$\frac16$ and so on - cannot be written as a decimal number [3].
In fact, only fraction whose denominator is of the form
$2^n\times5^m$, for some natural numbers $n$ and $m$, can be
written as decimal numbers. This may not be a serious problem but
it is a nuisance nonetheless.

If we want a nice decimal expression for $\frac13$ we must use a
system whose basis is a multiple of $3$, e.g. 3, 6 and so on.

A system with base $6$ (Fig. 2) would be very handy indeed. And
since this system has a very simple times table it is very
convenient too and would be welcomed in every school.

\bc
\begin{tabular}{c c c c c c c c c c c c c c c c }
  \hline &  ${}^|$&${}^|$ &  ${}^|$&${}^|$  & ${}^|$ & ${}^|$ & ${}^|$ &${}^|$
  &  ${}^|$ &${}^|$  &${}^|$  &${}^|$  &${}^|$  & ${}^|$&${}^|$  \\
 &$\>\>$ & $\>\>$ & $-\forall$ & $-$H
& $-\Pi$ & $-\Gamma$
& $\>0$ & $\>\Gamma$ & $\>\Pi$ & $\>\Delta$ & $\>$H
 & $\>\forall$ & $\Gamma0$&$\>\>$ &$\>\>$\\
& & & & & & & & & & & & & & &
\\
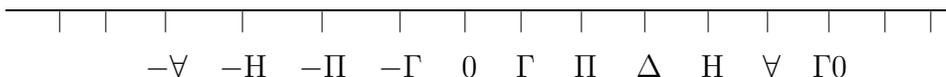
\multicolumn{14}{l}{\footnotesize{Figure 2: Number line for the
system
$N_{(0,5)}$}}\\
\end{tabular}
\ec

In fact, all the times tables in this system are either trivial -
the one for 0 and $\Gamma$ - or very easy - the one for $\Pi$,
$\Delta$ and $\forall$ (which play the role of the 9 in our number
system) with the exception of the table for H. There is therefore
a single difficult product, H$\times$ H. But we can easily move to
this point of the times table from other points, for instance, by
making H$\times$H = $\forall\times\Delta +\Gamma$. Therefore with
a system with basis 6, learning the times table would really be
child's play.

\section{Directed numerals and the joy of rounding}

Our figures can be misleading too. For example, 29 is far closer
to 30 than it is to 20 but its first numeral is 2. This problems
has psychological implication so that many shops use prices like
£199.99 because they surely believe (with good reason, no doubt)
that this figure conveys the impression of 'about one hundred
pounds' while it is obvious that the actual price is two hundred
pounds.

This problem is related to rounding. Rounding should be a natural
things to do with numbers. But we see that our number system is
not well behaved in this area.

Can this problem be solved? A simple solution - and one that, as a
byproduct, improves other areas of our number system too -
involves the introduction of the substraction operation in the way
we write numbers. The idea is to use numerals that instead of
adding, take away. Let, for instance, $\alpha$ be the symbol for
taking away one unit. Then we have $3\alpha = 30 - 1 = 29$.

Introducing these symbols may appear an artificial thing to do
since now we would have two different ways of writing the same
number. However, it is easy to see that we can write any number by
using just five negative numerals: $\alpha, \beta, \delta, \gamma$
and $\epsilon$, and four positive ones: 1, 2, 3, and 4 (Fig. 3).
(The negative numerals and the positive numerals can be
collectively referred to as 'directed numerals.')

 \bc
\begin{tabular}{c c c c c c c c c c c c c c c c }
  \hline &  ${}^|$&${}^|$ &  ${}^|$&${}^|$  & ${}^|$ & ${}^|$ & ${}^|$ &${}^|$
  &  ${}^|$ &${}^|$  &${}^|$  &${}^|$  &${}^|$  & ${}^|$&${}^|$  \\
 &$\alpha3$&$\alpha4$ & $\>\epsilon$ & $\> \gamma$ & $\>\delta$
& $\>\beta$ & $\>\alpha$
& $\>0$ & $\>1$ & $\>2$ & $\>3$ & $\>$4  & 1$\epsilon$ & 1$\gamma$&$1\delta$ \\
& & & & & & & & & & & & & & &
\\
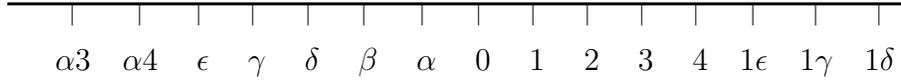
\multicolumn{14}{l}{\footnotesize{Figure 3: Number line for the
system
$N_{(5,4)}$}}\\
\end{tabular}
\ec

For example, we have now 284 =$3\beta4$. If we now round the last
numeral we get $3\beta0$. If we round now to hundreds we get 300.
So we see that this number system is better behaved  under
rounding than the current one.

Moreover, this number system brings with it unexpected gains, a
saving of symbols, for example. The sign '-' for indicating
negative numbers is no longer needed. With the new notation a
number will be negative if the first significant digit is negative
and positive if the first significant digit is positive. For
instance, $\delta3$ is negative ($=-27$).

A disadvantage of the new system with respect to the standard one
is its bad behavior under 'changing sign'. With the standard
number system this is done by simple adding (or removing) a '-'
sign. With the new system this is considerably more complicated.
For instance, the opposite of $2\epsilon$ is $\alpha\epsilon$.

This complication can be traced back to the lack of symmetry
between negative  numerals and positive ones. Whereas there are
five negative numerals there are only four positive ones [4]. This
problem would not exist if a system with an odd base were used.
For example with a base 7 made of three positive numerals, three
negative ones and zero.

\section{More generalizations}

The best way of understanding a structure is exploring ways in
which it could be changed while preserving its principal
characteristics. We have seen above that in addition to different
bases, number systems can be considered that make use of directed
numerals. This whole family of number systems can be characterized
by a pair of natural numbers $(a,b)$ where $a$ corresponds to the
number of negative numerals and $b$ to the number of positive
ones. In this way our number system would be $N_{(0,9)}$ and the
one described in the previous section would be $N_{(5,4)}$. A
symmetrical system - and  well behaved therefore under rounding -
with base 7 would be $N_{(3,3)}$.

Are there more radical generalization available? The following
definitions offer a wide family of possibilities.

\noindent {\it Definition 1.} A normal pre-numeration system is a
triplet $N=<I,\>i,\>(a)>$ where $I$ is a finite or countable set;
$i$ is an injective application from $I$ unto the set of integer
numbers $\cal{Z}$; and $(a)$ is a sequence valued in $\cal{Z}$.

The numerical value of a string (figure) $c_n...c_2c_1$ (with
$c_n,..., c_1 \in I )$ is given by

\be V_N(c_n...c_2c_1)=i(c_1)\times a_1 + i(c_2)\times a_2
+...+i(c_n)\times a_n \ee

It is clear that in this definition the elements of $I$ play the
role of the symbols (numerals) of the number system; the function
$i$ fixes the value of each numeral; and the sequence $(a)$ gives
a value to the numerals' place in the figure.

For our numeration system $I = \{ 0,1...9\}$; the application $i$
assigns the value 'nil' to 0, the value 'unity' to 1, and so on;
and the sequence $(a)$ is given by $a_n=10^n$.

\noindent{\it Definition 2.} A normal pre-numeration system
$N=<I,\>i,\>(a)>$ is said to be

\begin{itemize}

\item {\it Finite} iff $I$ is finite.
\item {\it Perfectly ordered} iff the sequence $(a)$ fulfil the
condition that $a_n = a_m \iff n=m $
\item {\it Perfectly disordered} iff the sequence $(a)$ is constant
[5].

\end{itemize}

Def. 1 allows us to easily consider 'exotic' number systems, where
positive and negative numerals are taken at random from the number
line. Other even more exotic possibilities - but which may be more
convenient to our psychology than the standard system -  include
using sequences such as

\be a_n = n! \qquad \hbox{so that, for instance, 23 means}\>
2\times 2! + 3\times 1 \ee or also

 \be a_n = n^n \qquad \hbox{so that 57 means}\>
5\times 2^2 + 7\times 1 \ee But, how many of these pre-numeration
systems are well behaved?

\noindent{\it Definition 3.} A normal pre-numeration system $N$ is
said to be a {\it normal number system} if every integer number
$n$ can be written within $N$, i.e. if for all integer $n$, there
exist a string of numerals $c_s...c_1$ such that $V_N(c_s...c_1) =
n$ [6].

\noindent{\it Definition 4.} A normal number system is said to be
{\it univocal} if every integer $n$ can be expressed in a single
way within $N$.

A exciting line of research would be to find out what kind of
normal pre-number systems are normal number systems and among
these which are univocal.

It is clear that generalizations that go beyond normal number
systems can also be considered.

There are three basic things that can be done with two numbers: 1)
compare them; 2) add them; and 3) multiply them. A number system
is a way of giving names to the numbers. The natural demand to be
made of a number system is that these names are well behaved under
these three operations. It turns out, however, that while our
number system is good with respect to the first operation
(comparing numbers) it is very badly behaved with respect to the
other two (adding and multiplying numbers). For there is nothing
in the pair of symbols 49, for instance, that indicates that this
number is $7\times7$. Or there is nothing in 18 that indicates
that it is worth $3\times6$. Or for that matter there is nothing
in 9 that indicates that it denotes 'nine'.

On the other hand, the most important property of a number system
is that the representation of a number conveys clearly the actual
size of the number. We have seem that not even in this area is our
current system specially good. On the contrary, it can be
misleading and is not well behaved under rounding.

A notation is good when it is used unnoticeably. This is clearly
not the case with our number system. Instead of facilitating our
calculations it is an obstacle that we have to climb to carry them
out.

From these consideration it follows that the shape of the
numerals(and the figures) deserves attention. Perhaps it would be
possible to devise numerals (and figures) that would multiply (and
add) naturally.

Ideally, it should be possible to graphically combine two figures
to obtain the figure representing their total or their product.

\section{Conclusions and perspectives}

We have discussed some of the positive aspects and some of the
negative aspects of our number system and have proposed ways of
overcoming the latter.

We have found that once generalizations are considered it becomes
apparent that our number system is not as magnificent at it
appears on the surface. Far from being perfect, our number system
has many problems and limitations.

An important problem of our number system is the complexity of its
times table. Many children start loathing Mathematics as soon as
they come across the times table. If the first thing our children
find of Mathematics is so unattractive, should we not take into
serious consideration the possibility of changing it?

We have suggested that a number system with base 6 is advantageous
in a number of ways. Its times table in particular is much easier
to learn. This system has the additional advantage that a system
in base 12 (dozens) is already used in some realms, e.g. in the
hours of the day, minutes in the hour, the months of the year, and
so on. Perhaps it would be advantageous to give more protagonism
to a system with base 6, which could be increasingly used in
everyday life and in the first years of Mathematics education.

We have also discussed other, more 'exotic' generalization of our
number system. Among these we have found several that solve some
of the problems of our number system. We have not found however
any which represent a significant improvement.

It is perfectly possible, however, that among the generalizations
yet to be analyzed a number system can be found that is really
advantageous. Perhaps one of these systems allows us to seamlessly
calculate 4,349 $\times$ 5,673 or $12.7\%$ of £15,642 and it is
perfect for everyday Mathematics. Perhaps another is most suitable
for advanced Mathematics, e.g. number theory.

It would not be inconceivable that alongside the current number
system another one (or several others) more appropriate for the
task at hand could be used.

It is clear that more research is needed in this area which has so
clear practical applications and so exciting theoretical
implications. The objective is devising a way of denoting the
numbers which is well behaved under addition and multiplication.
Ideally the figure denoting the product (or addition) of two
numbers should be obtainable by combining graphically (in an easy
way) the figures representing the factors. It is possible that
this can be done within the current decimal system by simply using
numerals with better, more suitable shapes.

Will our grand-grandchildren wonder in the future why despite all
our technological advances we have failed to perfect our
inconvenient number system for several centuries?

\section{Notes}

\noindent [1] In a sense a number system is just a way (more or
less convenient) of labelling (or naming) numbers. Thanks to the
place-value principle we can name $10^n$ numbers using figures
with no more than $n$ numerals. But numbers can also be denoted
with their name, in English, for instance.

\noindent [2] In fact, an interesting exercise would be to study
the possibility (theoretical at least) of constructing a
perfectly-ordered language, i.e. a language such that the order of
the words plays a fundamental role in conveying meaning and such
that a meaning could be conveyed with one order only.

\noindent [3] Sometimes an arc over the recurring digits is used.
But in purity this, in addition to being inelegant, implies
actually going beyond the place value principle.

\noindent [4] This asymmetry implies that the new number system is
not completely well behaved under rounding.  For instance,
$2\epsilon\epsilon\> (=145)$ is rounded (hundreds) to 100 despite
its first digit being 2.

\noindent [5] Originally the order of the symbols was completely
irrelevant in the Roman number system and $VI$, for instance, and
$IV$ denoted the same number, namely 'six'. The original Roman
number system was therefore completely disordered. The (minor)
role of the order was introduced in the Middle Ages (and it is not
clear at all that it represents an improvement).

\noindent [6] A subtlety may complicate this definition. In fact
it is only necessary that the figures cover half the integers - as
long as the other half can be reached by using a negative '-'
sign, for instance.
    In this definition we have not paid attention either to the
decimal numbers.

\noindent [7] A base 6 appears natural too from a Mathematical
point of view in that 6 is the product of the first two prime
numbers. According to this criterion the possible bases would be
2, $6=2\times3$, $30 = 2\times3\times5$, 210... and so on.

\end{document}